\documentstyle[12pt,amstex,amscd,amssymb]{article} 
\newtheorem{theorem}{Theorem}
\newtheorem{lemma}{Lemma}

\begin{document}

 \title
{Threefolds with Vanishing Hodge Cohomology}  
\author
{Jing  Zhang}
\date{}
\maketitle
\begin{abstract}
We consider algebraic manifolds $Y$ of dimension 3 
over  $\Bbb{C}$ with $H^i(Y, \Omega^j_Y)=0$
for all $j\geq 0$ and $i>0$. Let $X$ be a smooth completion
of $Y$ with
 $D=X-Y$, an effective divisor on $X$
with  normal crossings. If the $D$-dimension of $X$ is not zero,
then $Y$ is a fibre space over a smooth affine curve $C$ (i.e.,  we 
have a surjective morphism from 
$Y$ to $C$ such that general fibre is smooth and irreducible) 
such that   every fibre satisfies the same vanishing condition.  
If an irreducible smooth fibre is not affine, 
 then the Kodaira dimension of $X$
is $-\infty$ and the $D$-dimension of X is  1. We also discuss 
sufficient conditions from the behavior of fibres or higher direct 
images to guarantee the global vanishing of Hodge cohomology and the 
affineness of $Y$.    
\end{abstract}

AMS(1991)Classification. 14J30, 14B15, 14C20.

 \begin{center}
\large {0. Introduction}
\end{center}

Let $Y$ be a complex manifold
with   $H^i(Y, \Omega^j_Y)=0$   for all $j\geq 0$ and $i>0$, 
then what is $Y$? Is Y Stein? This is a question raised by Serre [Se]. 
 Peternell [P] also asked the same question for schemes: 
 If $Y$ is  a smooth   scheme
of finite type over $\Bbb{C}$, 
 is it affine?  For the  nonalgebraic case, in particular, 
 complex surfaces,  see  Peternell's paper [P]. 
 Throughout this paper, we assume that $Y$ is an algebraic manifold,
i.e., an irreducible nonsingular algebraic variety defined over $\Bbb{C}$
. If dim$Y=1$, then $Y$
is affine. If dim$Y=2$, Mohan Kumar [Ku] classified it completely. 
It may not be 
affine and  has three possibilities as follows:

      (1) $Y$ is affine.

      (2) Let $C$ be an elliptic curve and $E$ the unique nonsplit 
      extension of $\mathcal{O}$$_C$ by itself.  
      Let ${X=\Bbb{P}}_C(E)$ and  $D$ be the canonical section, then $Y=X-D$.

       (3) Let $X$ be a projective rational surface with an effective 
       divisor $D=-K$ with $D^2=0$, $\mathcal{O}$$(D)|_D$ be nontorsion and 
       the dual graph of $D$ be $\tilde{D}_8$ or $\tilde{E}_8$, then $Y=X-D$. 

 If the surface Y is not affine, then the Kodaira dimension of
 X is  $-\infty$, $D-$dimension is 0 ([Ku], Lemma 1.8). 
 In the second case, the 
canonical divisor $K_X=-2D$ and $D$ is irreducible, so the  logarithmic Kodaira
 dimension  $\bar{\kappa}(Y)=\kappa(D+K_X, X)=-\infty$; in the third case, 
 $K_X=-D$ and $D$ is either $\tilde{D}_8$ or $\tilde{E}_8$.  
 Let $D'=\sum D_i$ be the reduced divisor, where $D_i's$ are the
 prime  components of $D$,
 then $K_X+D'$ is not effective. Therefore again 
   $\bar{\kappa}(Y)=-\infty$ ([I3], [Ku], [Mi]).
  The surfaces in the first two cases are Stein. 
  The third case is open.   
  Since there exist Stein varieties which are not affine 
  (see [H2], Serre gave the first example of this),  
  the surfaces in the third case might be Stein.

       If dim$Y=3$, let us first fix our basic assumption (BA) as follows.

{\bf{(BA)}} Let Y be a smooth  irreducible threefold,  $X$ be
 a smooth completion (See Nagata [N] for the existence of $X$)
  such that  $X-Y=\cup D_i$, union of connected, distinct 
  prime divisors on $X$. Let $D$ be an effective divisor supported 
  in $\cup D_i$ with  normal crossings [I3]. Suppose that 
  the moving part of $|nD|$ is base point free (such a smooth
  completion of $Y$ always exists after further blowups),  
  $H^0(X, \mathcal{O}$$_X(nD))\not=\Bbb{C}$ for some $n$ and $Y$ 
  contains no complete surfaces.

  {\bf{Proposition}}
{\it  Under the condition (BA), there is a smooth projective 
curve $\bar{C}$,  and a smooth  affine curve $C$ such that 
the following diagram commutes

\[
  \begin{array}{ccc}
    Y                           &
     {\hookrightarrow} &
    X                                 \\
    \Big\downarrow\vcenter{%
        \rlap{$\scriptstyle{f|_Y}$}}              &  &
    \Big\downarrow\vcenter{%
       \rlap{$\scriptstyle{f}$}}      \\
C        & \hookrightarrow &
\bar{C}
\end{array}
\]
where f is proper and surjective, every fibre of f over $\bar{C}$
 is connected, 
 general fibre  is  smooth. Also general fibre 
 of $f|_Y$ is connected  and smooth. And if 
 the $D-$dimension of $X$
is no less than 2, then we can make  $\bar{C}$ to be ${\Bbb{P}}^1$.}

Our main results are the following.

{\bf{Theorem A}}
{\it If  $H^i(Y, \Omega^j_Y)=0$ for all $j\geq 0$ and $i>0$ and 
$H^0(Y, \mathcal{O}$$_Y)\not=\Bbb{C}$, then for the above $f$, we have

       (1)Every fibre $S$ of $f|_Y$ over $C$  satisfies the same vanishing 
       condition, i.e., 
$H^i(S, \Omega^j_S)=0$.

         (2) If there is a smooth fibre $X_{y_0}$ over  
 $y_0\in \bar{C}$ such that  $X_{y_0}|_Y=S_0$
is not affine, then the Kodaira dimension of $X$ is $-\infty$ and 
the $D$-dimension of $X$ is 1.

(3) If one fibre $S$ of $f|_Y$ over $C$  is not affine, 
then $Y$ is not affine.
$Y$ is affine if and only if for every coherent sheaf $F$ on $X$
$$  h^1(X,\lim_{{\stackrel{\to}{n}}}
F\otimes {\mathcal{O}}(nD)) <\infty.
$$}

        Conversely, in the above diagram, 
        let $F_n=\Omega^j_X\otimes \mathcal{O}$$(nD)$ or
         $F_n=\Omega^j_X$(log$D)\otimes \mathcal{O}$$(nD)$,
where $\Omega^j_X$(log$D)$ is the sheaf of logarithmic $j$-forms on 
$Y$ ([I1], [I2], [I3]), then we have

{\bf{Theorem B}} {\it  If  the higher direct images satisfy  
$$     \lim_{{\stackrel{\to}{n}}}
      R^1f_*F_n=
\lim_{{\stackrel{\to}{n}}}R^2f_*F_n=0,$$
or for every point $y\in C$, if  $D_y=X_y\cap D$ is a curve, and 
  $$\lim_{{\stackrel{\to}{n}}} H^2(X_y, F_{n,y})=0, \quad 
  \quad \lim_{{\stackrel{\to}{n}}}R^1f_*F_n=0,
$$
where $F_{n,y}=F_n|_{X_y},$ 
 then $H^i(Y, \Omega^j_Y)=0$ for all $j\geq 0$ and $i>0$.}

   Similar results for the affineness of $Y$ can be obtained. We will 
   discuss it in section 3. 
   From theorem $A$ and the surface cases (2) and (3)
which are not affine, we know that a threefold with the vanishing 
Hodge cohomology
is not necessarily affine. If it is affine, then of course it is Stein.
In surface case, if $Y$ is not affine, the Kodaira dimension of $X$ 
is unique. 
Is it still true for threefolds? When there are nonconstant regular functions 
on the threefold $Y$, by theorem $A$(1) and Mohan Kumar's classification, 
there are three 
different types of smooth fibres.   
The question reduces to two questions. 
The first question is: given a smooth variety $Y$ 
with $H^i(Y, \Omega^j_Y)=0$
for all $j\geq 0$, $i>0$ and a surjective morphism from $Y$ to a smooth 
affine curve $C$ such that 
every (or general) fibre is affine, then is $Y$ 
affine? Generally it is not true (even in the surface case) without restriction 
on Hodge cohomology. Under our cohomology restrictions, if $Y$ is a surface, 
then it is true (Lemma 1.8, [Ku]). The second question 
is the invariance of plurigenera. 
If one fibre $S_0$ is not affine, then is the Kodaira dimension 
of $X_t$ a constant 
in an open neighborhood of $0$? Iitaka conjectured that in a
 smooth family the $m$th
plurigenus is constant. He proved it for surface case ([I4], [I5]).
 Nakayama 
proved that the conjecture follows from the minimal model
 conjecture and the 
abundance conjecture ([Na1], [Na3]). Siu proved it if the
generic fibre is of 
 general type [Si].
Kawamata extended Siu's result to fibres with canonical singularities
[Ka4].   In our case, some isolated fibres may be singular or reducible 
or both. We can not therefore apply these results.

Now let $Y$ and $X$  be as in theorem $A$. We are sure that the Kodaira 
dimension of $X$
can be $-\infty$  and the $D$-dimension can be 1 (Theorem 7).
 Our motivation 
here is  to see the global picture from the fibre, i.e., if every fibre 
or general fibre has vanishing Hodge cohomology, then is it still true 
for $Y$?
We can prove that the direct limit of second direct images vanishes
thus  $H^2(Y, \Omega^j_Y)=0$. But
  the direct limit of the first direct images might be supported at 
  finitely many points. We do not know how to deal with these points 
  and whether the first  direct image sheaves are locally 
  free or not. 
In fact, by a  result of Goodman and Hartshorne (Lemma 4),  we only 
need the local freeness on $C$.
If it is true, then the direct limit is zero and therefore we also have 
$H^1(Y, \Omega^j_Y)=0$. Hence we can get an equivalent condition for 
vanishing 
Hodge cohomology of $Y$.

 Mohan Kumar's proof in surface case heavily depends on 
 the following two facts.
 By his Lemma 1.10, any line bundle $L$ on the  carefully chosen divisor $D$   
  with degree zero when restricted to each component of $D$ 
has the following  property: $H^0(L)\neq 0$ if and only if $L\cong 
{\mathcal{O}}_D.$ Thus if $L={\mathcal{O}}_D(D)$, then it satisfies all
 these conditions by the choice of $D$ therefore it  is either torsion 
or nontorsion. This is why $Y$ has only two possibilities if it is not
 affine.
We have no any idea of a similar result when the dimension is 3.  
The second fact is   Zariski decomposition of $D$. He used it
 to compute intersection numbers and $h^0(X, {\mathcal{O}}_X(mD))=1$
for every nonnegative $m$. 
 But in threefold, we do not always have  Zariski decomposition
 (see [C]. 
For recent progress, see [Na2]).  To understand $Y$ and $X$, 
we have to use a different approach. 
 We first  construct a proper and surjective morphism from $X$ to a 
 smooth curve. This is the place
where  we need the condition that $D$-dimension of $X$ is not 0. 
Notice that we can not use any other divisor on $X$ to define our map. 
Otherwise, we have no control of the cohomology and 
the image of $Y$ then can not use our assumption. 
This is saying that 
we can only change the boundary $D$ but can not change $Y$.
This is why we can not use Iitaka's fibration and Mori's construction. 
And in order to use Iitaka's fibration, we must assume 
$\kappa (X)\geq 0$ which is not true in our case.
 But we can use Iitaka's $C_n$ conjecture ([Ka2], [V]).

        The content of this paper is divided into three parts. 
        In the first section, we will present some basic lemmas. 
        We borrow the idea from the surface case, i.e., from  [Ku], [P].
 The second section contains construction of the fibre space. We will
  prove our main theorems in the third section and give an example. 
  Our basic tools are Grothendick's local cohomology theory 
  and classification theory developed by the Japanese school of algebraic 
  geometry. 

     {\bf  Convention} Unless otherwise explicitly mentioned, we always 
     use Zariski topology, i.e.,  an open set means a Zariski open set.

{\bf{Acknowledgments}} \quad \quad It is my great pleasure to thank 
 my adviser, Professor N. Mohan Kumar
to lead me to the beautiful field of algebraic geometry, in particular, 
to the fascinating theory of threefolds.
 I have benefited from  communication with the following professors: 
 Michael Artin,  Quo-Shin Chi, Steven D. Cutkosky, Xiaojun Huang, 
  Nicholas M. Katz, Yujiro Kawamata, 
 J$\mbox{\'{a}}$nos Koll$\mbox{\'{a}}$r,
Tie Luo, Noboru Nakayama,  Zhenbo Qin, A. Prabhakar Rao,  
 David Wright  and   Qi Zhang. 
I would like to thank all of them.  I also thank the referee for his/her
comments.

 \begin{center}
\large 1. Preliminary lemmas
\end{center}

\begin{lemma} Let $Y$ be an irreducible smooth  threefold with 
 $H^i(Y,\Omega^j_Y)=0$
for every $j\geq 0$ and $i>0$. Let 
$X$ be any smooth 
completion of $Y$, then
$X-Y$ has no isolated points. 
\end{lemma}
$Proof$.  If $P$ is an isolated point of $X-Y$, let $Y'=Y\cup \{P\},$ then $Y'$
is a scheme and we have  exact sequence of local cohomology  
$$ 
0=
H^2(Y, {\mathcal{O}}_Y)
\longrightarrow 
H^3_{\{P\}}({\mathcal{O}}_{Y'})
\longrightarrow 
H^3(Y', {\mathcal{O}}_{Y'})
=0,
$$
The last term is zero since $Y'$ is not complete.
But 
$$
H^3_{\{P\}}({\mathcal{O}}_{Y'})
\cong 
\lim_{{\stackrel{\to}{n}}}
Ext^3_{{\mathcal{O}}_{X}}(
{\mathcal{O}}_{nP}, {\mathcal{O}}_{Y'})
\neq
0
$$
where ${\mathcal{O}}_{nP}
={\mathcal{O}}_{X}/{\mathcal{M}}^n$,
$\mathcal{M}$ is the ideal sheaf of $P$. 
To see this, write the  short exact sequence
$$0
\longrightarrow
{\mathcal{M}}^n/{\mathcal{M}}^{n+1}
\longrightarrow
{\mathcal{O}}_{(n+1)P}
\longrightarrow
{\mathcal{O}}_{nP}
\longrightarrow
0.
$$
Then,  since $Ext^2_{{\mathcal{O}}_{X}}
({\mathcal{M}}^n/{\mathcal{M}}^{n+1}
, {\mathcal{O}}_{Y'})=0, 
$ we have 
$$0
\longrightarrow
Ext^3_{{\mathcal{O}}_{X}}(
{\mathcal{O}}_{nP}, {\mathcal{O}}_{Y'})
\longrightarrow
Ext^3_{{\mathcal{O}}_{X}}(
{\mathcal{O}}_{(n+1)P}, {\mathcal{O}}_{Y'})
\longrightarrow
Ext^3_{{\mathcal{O}}_{X}}
({\mathcal{M}}^n/{\mathcal{M}}^{n+1}
, {\mathcal{O}}_{Y'})
\longrightarrow
0.
$$
If $Ext^3_{{\mathcal{O}}_{X}}
({\mathcal{M}}^n/{\mathcal{M}}^{n+1}
, {\mathcal{O}}_{Y'})
\neq 0$,
then dimension of $
Ext^3_{{\mathcal{O}}_{X}}(
{\mathcal{O}}_{nP}, {\mathcal{O}}_{Y'})
\longrightarrow  \infty$ as $n\longrightarrow \infty$.
Thus  $H^3_{\{P\}}({\mathcal{O}}_{Y'})
\neq 0$. For some suitable $m$
determined by $n$, we have 
$${\mathcal{M}}^n/{\mathcal{M}}^{n+1}
=({\mathcal{O}}_{P}/{\mathcal{M}})^m.
$$ 
Therefore 
$$
Ext^3_{{\mathcal{O}}_{X}}
({\mathcal{M}}^n/{\mathcal{M}}^{n+1}
, {\mathcal{O}}_{Y'})
=\oplus Ext^3({\mathcal{O}}_{P}/{\mathcal{M}}, 
{\mathcal{O}}_{Y'})
= \oplus Ext^3({\Bbb{C}}(P),{\mathcal{O}}_{Y'}).
$$ 
Choose local coordinates such that $P=\{x=y=z=0\}, x, y, z\in 
{\mathcal{O}}_{U}, U$ is a neighborhood of $P$, then 
$$0
\longrightarrow
{\mathcal{O}}_{U}
\longrightarrow
{\mathcal{O}}_{U}^3
\longrightarrow
{\mathcal{O}}_{U}^3
\longrightarrow
{\mathcal{O}}_{U}
\longrightarrow
{\Bbb{C}}(P)
\longrightarrow
0.
$$
Moreover, 
\[{\mathcal{E}}xt^i({\Bbb{C}}(P), {\mathcal{O}}_U)
=\left\{ \begin{array}{ll}
0            & \mbox{if $i\neq 3$}  \\
{\Bbb{C}}(P) &  \mbox{i=3}.
   \end{array}
\right.   \]
Finally we can compute 
$$
Ext^3({\Bbb{C}}(P),{\mathcal{O}}_{Y'})
=H^0(
{\mathcal{E}}xt^3
({\Bbb{C}}(P), {\mathcal{O}}_U))
={\Bbb{C}}(P) \neq 0.
$$
  \begin{flushright}
 Q.E.D. 
\end{flushright}

\begin{lemma} Under the condition of Lemma 1, $X-Y=\cup Z_i$ is connected, 
where  $Z_i$'s are irreducible  components.
\end{lemma}
$Proof$.  If $X-Y=Z$ is not connected, write $Z=Z_1+Z_2$,  $Z_1\cap    
     Z_2=\emptyset.$ We have a long exact sequence of local cohomology
$$0=H^2(Y,\Omega^3_Y)\longrightarrow 
      H^3_Z(X,\Omega^3_X)\longrightarrow 
 H^3(X,\Omega^3_X)\longrightarrow  H^3(Y,\Omega^3_Y)
=0.$$
So $ H^3_Z(X,\Omega^3_X)=H^3(X,\Omega^3_X)=\Bbb{C}$ by Serre duality. 
But by Mayer-Vietoris sequence,
$$ H^3_Z(X,\Omega^3_X)\cong  H^3_{Z_1}(X,\Omega^3_X)\oplus 
 H^3_{Z_2}(X,\Omega^3_X).$$ 
Both summands are at least one dimensional since 
$$ H^3_{Z_i}(X,\Omega^3_X)\longrightarrow  H^3(X,\Omega^3_X)={\Bbb{C}}
 \longrightarrow    H^3(X-Z_i,\Omega^3_X)=0.$$
This is a contradiction.   
\begin{flushright}
 Q.E.D. 
\end{flushright}

 \begin{lemma} Let $X$, $Y$ be as above, then $Y$ contains no complete surfaces. 
 \end{lemma}
       $Proof$. If $S$ is a complete, irreducible surface in $Y$, 
       then we have short exact sequence
$$ 0\longrightarrow A\longrightarrow \Omega^2_Y 
 \longrightarrow \Omega^2_S \longrightarrow 0$$
where $A$ is the kernel. Since 
$\Omega^2_Y$ and  $\Omega^2_S$
are coherent, $A$ is coherent. 
For any abelian sheaf $\mathcal{F}$ on $X$, we have long exact sequence
$$...\longrightarrow  H^3_Z(X, {\mathcal{F}})  \longrightarrow 
   H^3(X, {\mathcal{F}})      \longrightarrow H^3(Y, {\mathcal{F}})\longrightarrow  0. $$
By formal duality [H3], 
$$H^0(\hat{X}, \hat{\mathcal{G}})=H^3_Z(X,\mathcal{H})^*$$
where ${\mathcal{G}}=
{{\mathcal{H}}{om}}_{{\mathcal{O}}_X}({\mathcal{F}}, \omega)$,
$\omega=\Omega^3_X$ and 
${\mathcal{H}}=
{{\mathcal{H}}om}_{{\mathcal{O}}_X}
(\mathcal{G}, \omega)$.
 If  $\mathcal{F}$ is locally free, then 
 ${\mathcal{H}}=F$. So we have
$$H^0(\hat{X}, \hat{\mathcal{G}})=H^3_Z(X,\mathcal{F})^*.$$
But   $H^0(X, {\mathcal{G}}) \longrightarrow  
H^0(\hat{X}, \hat{\mathcal{G}})$ 
is injective, by Serre duality, 
$H^3_Z(X, {\mathcal{F}})
\longrightarrow 
H^3(X, \mathcal{F})$
is surjective. So 
$H^3(Y, {\mathcal{F}})=0$
for any locally free sheaf $\mathcal{F}$. 
Then for any coherent sheaf  $\mathcal{F}$, 
$H^3(Y, {\mathcal{F}})=0$
since we have short exact sequence
$$0\longrightarrow B \longrightarrow {\mathcal{F'}}  
\longrightarrow  {\mathcal{F}}\longrightarrow  0$$
where   $\mathcal{F'}$ is locally free. 
In particular, $H^3(Y,A)=0$.  
From
$$0=H^2(Y,\Omega^2_Y)
\longrightarrow 
H^2(\Omega^2_S) \longrightarrow 
H^3(Y,A)=0$$ 
we have  $H^2(\Omega^2_S)=0$, which is a contradiction ([AK]). 
\begin{flushright}
 Q.E.D. 
\end{flushright}

{\bf{Remark 1}} Our proof is algebraic. In analytic category, 
we can use Siu's theorem [Si1] to get $H^3(Y, F)=0$ for any
analytic sheaf $F$ since $Y$ is not compact by Serre duality.
Then we can use
 Norguet and Siu's result ([NS], [P]). It says that if a complex 
manifold $Y$ contains a compact analytic subvariety of dimension
$q$ and for every coherent sheaf $F$ on $Y$, $H^{q+1}(Y, F)=0$,
then  $H^q(Y, \Omega^q)\neq 0$.

By the above lemmas, we know that for any smooth completion $X$ of $Y$,
the dimension of the boundary $X-Y$ is not zero. By suitable blowing ups,
we may assume they satisfy the basic assumption  (BA). So the $D-$dimension 
of such  $X$ 
makes sense. We put a lemma after theorem 1 for logical correctness. It says 
that if the $D$-dimension of $X$ is not zero, then $Y$ contains no complete
curves.

\begin{lemma}[Goodman, Hartshorne]  Let $V$ be a scheme and  
$D$ be an effective Cartier divisor on  $V$. Let $U=V-$Supp$D$ and $F$ 
be any coherent sheaf on $V$, 
then for every $i\geq 0,$ 
$$\lim_{{\stackrel{\to}{n}}}
H^i(V, F\otimes {\mathcal{O}}(nD)) \cong  H^i(U,  F|_U).
$$
\end{lemma}

 \begin{center}
\large 2. Construction of a proper, surjective morphism from $X$ to $\bar{C}$
\end{center}
  
If $H^0(X, {\mathcal{O}}_X(nD))\neq \Bbb{C},$ let $\xi$ be a nonconstant, 
irreducible element 
(which means that it can not be written as a product of two 
nonconstant elements)
 in  $H^0(X, {\mathcal{O}}_X(nD)),$ then it defines 
a rational map  
$$   \xi:  X\cdot\cdot \longrightarrow {\Bbb{P}}^1,
$$
with poles in $D$. When restricted to $Y$, it is a morphism
$$   \xi |_Y:     Y\longrightarrow {\Bbb{A}}^1.
$$
Let $U=\xi (Y)$, the image of $Y$ under $\xi$.
 By Hironaka's elimination of indeterminacy,  
 there is a smooth  projective variety $\tilde{X}$, such 
 that the morphism $\sigma : \tilde{X}\rightarrow X$ is 
 composite of finitely many monoidal transformations
which is isomorphic when restricted to $Y$, i.e., $Y$ is fixed  and 
$g=\xi\circ \sigma: \tilde{X}\rightarrow {\Bbb{P}}^1$ is proper
and  surjective. Replace $X$ by $\tilde{X}$, 
since $g|_Y=\xi,$ we have a commutative diagram
\[
  \begin{array}{ccc}
    Y                           &
     {\hookrightarrow} &
    X                                 \\
    \Big\downarrow\vcenter{%
        \rlap{$\scriptstyle{g|_Y}$}}              &  &
    \Big\downarrow\vcenter{%
       \rlap{$\scriptstyle{g}$}}      \\
U        & \hookrightarrow &
 {\Bbb{P}}^1.
\end{array}
\]
 To guarantee
the connectedness of fibres, we can use Stein factorization. Let 
$f: X\rightarrow \bar{C}$ 
be a proper surjective morphism, 
$h: \bar{C}\rightarrow {\Bbb{P}}^1$ 
be finite ramified covering  such that 
$g$ is composition of these two maps, i.e., $g=h\circ f$. 
Let $C=f(Y)$, then we have commutative diagram
\[
  \begin{array}{ccc}
    Y                           &
     {\hookrightarrow} &
    X                                 \\
    \Big\downarrow\vcenter{%
        \rlap{$\scriptstyle{f|_Y}$}}              &  &
    \Big\downarrow\vcenter{%
       \rlap{$\scriptstyle{f}$}}      \\
C        & \hookrightarrow &
\bar{C},
\end{array}
\]
where $f$ is proper and surjective and every fibre of $f$ 
is connected. Moreover,  $C$ and $\bar{C}$ are smooth.

Now consider the image of $D$ under $f$. If $f(D)$ is a point, 
then $Y$ contains  complete  surfaces, so $f(D)=\bar{C}.$ 
Since both $X$ and  $\bar{C}$ are irreducible and projective, every fibre of 
$f$ over  $\bar{C}$ has dimension at least 2 ([Sh], Chapter 1, section 6.3,
Theorem 7) but can not be 3 ([Sh], Chapter 1, section 6.1, Theorem 1),
that is, every fibre has dimension 2. 
By the second Bertini theorem([Sh], Chapter 2, section 6.2), 
there is an open set $U\subset \bar{C}$, 
such that every fibre $f^{-1}(P)$ for every point
$P$ in $U$ is smooth.

Since $f(D)=\bar{C},$  there is a component $D_i$ of $D$, 
such that $f(D_i)=\bar{C}.$
But some components of $D$ may have points as images. 
Removing these finitely many points from $C,$ 
for general point $P$ in $C$, the inverse image $\bar{S}=f^{-1}(P)$
is an irreducible surface such that $\emptyset\neq D_i\cap \bar{S}\subset
 D\cap \bar{S}.$ By irreducibility of general fibre,  $ D\cap \bar{S}$ is
  a curve on $\bar{S}$ for general 
$P$. Removing this 
curve, the surface $S=\bar{S}-D=\bar{S}\cap Y$ is irreducible. So general 
fibre of $f|_Y$ over $C$ is smooth and irreducible(thus connected). 

By our construction,  $f_*{\mathcal{O}}_X={\mathcal{O}}_{\bar{C}}$
([U2], Proposition 1.13). But we do not know what the curve 
$\bar{C}$ is. If the $D$-dimension $\kappa (D)\geq 2,$ then we can make 
$\bar{C}$ to be ${\Bbb{P}}^1.$  The construction of rational map 
from $X$ to ${\Bbb{P}}^1$ is due to Ueno ([U2], page 46).

Choose two algebraically independent rational functions 
$\eta_1$ and $\eta_2$
in ${\Bbb{C}}(X)$. 
By Zariski's lemma ([HP], Chapter X, section 13, Theorem 1, page 78), 
there exists a constant $d$ such that the field 
${\Bbb{C}}(\eta_1 +d\eta_2)$ is algebraically closed in ${\Bbb{C}}(X)$. 
Define a rational map $f$ from $X$ to  ${\Bbb{P}}^1$ by sending points 
$x$ in $X$ to 
$(1, \eta_1(x) +d\eta_2(x) )$ in  ${\Bbb{P}}^1.$  
We can choose 
 $\eta_1$ and $\eta_2$, 
such that  
$\eta_1 +d\eta_2$ 
only has poles in $D$ ([U2], Lemma 4.20.3), that is, 
when restricted to $Y$, $f$ is morphism, 
then by our previous argument, we have diagram 
\[
  \begin{array}{ccc}
    Y                           &
     {\hookrightarrow} &
    X                                 \\
    \Big\downarrow\vcenter{%
        \rlap{$\scriptstyle{f|_Y}$}}              &  &
    \Big\downarrow\vcenter{%
       \rlap{$\scriptstyle{f}$}}      \\
C        & \hookrightarrow &
{\Bbb{P}}^1,
\end{array}
\]
where $f$ and $f|_Y$ satisfy the same properties as before.

{\bf{Proposition}}
{\it       Under the condition (BA), there is a smooth 
projective curve $\bar{C}$, 
        and a smooth, affine curve $C$ such that the following diagram
         commutes

\[
  \begin{array}{ccc}
    Y                           &
     {\hookrightarrow} &
    X                                 \\
    \Big\downarrow\vcenter{%
        \rlap{$\scriptstyle{f|_Y}$}}              &  &
    \Big\downarrow\vcenter{%
       \rlap{$\scriptstyle{f}$}}      \\
C        & \hookrightarrow &
\bar{C}
\end{array}
\]
where $f$ is proper and surjective, every fibre of $f$ over $\bar{C}$ is 
connected, general fibre of $f$ is  smooth. Also general fibre of $f|_Y$ 
is connected and smooth. Moreover, if the $D$-dimension 
of $X$ is no less than 2, then we can make $\bar{C}$ to be  ${\Bbb{P}}^1$. 
}

{\bf{Remark 2}} By [I1], page 79, for general fibre $X_y=f^{-1}(y)$, 
$D|_{X_y}=D_y$ is a divisor of $X_y$ with normal crossings 
if $D$ is a divisor 
with normal crossings.

 \begin{center}
\large 3. Structure of $Y$ with $h^0(X, \mathcal{O}$$_X(nD))>1$ 
\end{center}

In the diagram of the proposition, since $\bar{C}$ is smooth, $f$ is flat.

 \begin{theorem} If $Y$ is a smooth three-fold with 
 $H^i(Y, \Omega^j_Y)=0$   for all $j\geq 0$ and $i>0$ 
and  $H^0(X, \mathcal{O}$$_X(nD))\not=\Bbb{C}$ for some $n$,
then in the construction of the proposition, for every fibre
 $S$ of $f|_Y$ over $C$, 
$H^i(S, \Omega^j_Y|_S)=0$   therefore  
$H^i(S, \Omega^j_S)=0$   for all $j\geq 0$ and $i>0$. 
\end{theorem}
$Proof$. $Y$ has no complete surfaces by lemma 1. Thus the
 condition of the proposition is satisfied. For any point $P$ on $C$, let 
$g$ be an  element of $\Gamma (C, {\mathcal{O}}_C)$ such that the divisor 
defined by $g$ is 
$Q$=div$g=P+Q_1+\cdot\cdot\cdot +Q_r$, $P\neq Q_i$ for every $i,$ then 
$$    f^{-1}(Q)=S_Q=S\cup S_1\cup\cdot\cdot\cdot\cup S_r
$$   
where $S$ is the fibre over $P$, $S_i$ is the fibre over $Q_i$. 
From the short exact sequence
$$ 0\longrightarrow 
 {\mathcal{O}}_Y
\longrightarrow 
 {\mathcal{O}}_Y
\longrightarrow 
 {\mathcal{O}}_{S_Q}
\longrightarrow 
0
$$
where the first map is defined by $g$, we have 
$$     H^i(S_Q,   {\mathcal{O}}_{S_Q})=0
$$
for every $i>0$. 
Similarly, from the short exact sequence  
$$ 0\longrightarrow 
 \Omega^j_Y
\longrightarrow 
 \Omega^j_Y
\longrightarrow 
 {\Omega^j_Y}|_{S_Q}
\longrightarrow 
0
$$
where the first map is still defined by $g$, 
we have 
$$    H^i(S_Q,  {\Omega^j_Y}|_{S_Q})=0.
$$
By Mayer-Vietoris sequence, we have 
$$
H^i(S, \Omega^j_Y|_S)=0.
$$
In particular, $H^i(S, {\mathcal{O}}_S)=0$. 
 From the  exact sequence  
$$0\longrightarrow 
A\longrightarrow
  \Omega^j_Y|_S
\longrightarrow 
\Omega^j_S
\longrightarrow 
0
$$
 we have $H^i(S, \Omega^j_S)=0$ for every $i>0$
and  $j\geq 0$ since 
 $H^2(S, A)=0$ for coherent sheaf $A$ ([H3], [Kl]).  
             
 \begin{flushright}
 Q.E.D. 
\end{flushright}

{\bf{Remark 3}} If $H^i(Y, \Omega^j_Y)=0$ for every $i>0$ and $j\geq 0$, 
$H^0(X, \mathcal{O}$$_X(nD))\neq 
\Bbb{C}$ for some $n$, then $Y$ contains no complete curves. In fact, 
if $E$ is such a curve 
in $Y$, then its image under $f|_Y$ is a point $P$ on $C$, 
so $E$ is contained in the fibre $S$ of $f|_Y$ over $P$ also contained 
in $f^{-1}(P)=X_P$ in $X$. Write $X_P=X_P'+D'$ where $D'$ is a divisor 
contained in $D$
and  $X_P'$ intersects $Y$ with the surface $S$ in $Y$, i.e., 
$S=Y\cap X_P'$, 
and  $X_P'\cap D$ is a curve, 
then  $H^i(S, \Omega^j_S)=0$ for every $i>0$ and $j\geq 0$ by theorem 1. 
 This implies that $S$ is not complete [AK]. 
If there is a complete curve $Z$ in $S$, then $H^1(Z, \Omega^1_Z)\neq 0$  but
$$ 0\longrightarrow A\longrightarrow \Omega^1_S 
 \longrightarrow \Omega^1_Z\longrightarrow 0$$
and $H^2(S, A)=0$ ([H3], [Kl]). This is a contradiction. So we have showed

\begin{lemma} Let $Y$ be a smooth threefold with  $H^i(Y, \Omega^j_Y)=0$ 
for all $j\geq 0$ and $i>0$ and 
$H^0(Y, \mathcal{O}$$_Y)\not=\Bbb{C}$, then $Y$ contains no complete curves.
\end{lemma}

Now consider the sheaves 
$ \Omega^j_X$ and 
$\Omega^j_X$(log$D)$. 
Let $F_n=\Omega^j_X\otimes \mathcal{O}$$(nD)$ or 
$F_n=\Omega^j_X$(log$D)\otimes \mathcal{O}$$(nD)$, 
then $F_n$  is flat over $\bar{C}$ since it is locally free on $X$
 and $\bar{C}$ is smooth.  If $H^i(Y, \Omega^j_Y)=0$, then for every fibre
  $S=f_{|Y}^{-1}(y)=X_y\cap Y$, $y\in C$,  
 $H^i(S, \Omega^j_Y|_S)=0$  
for all $j\geq 0$ and $i>0$.  If $X_y$ is irreducible, then by lemma 4,
since  $F_n|_S=\Omega^j_Y|_S$, $D|_{X_y}=D_y$ is a divisor on $X_y$, we have
$$ \lim_{{\stackrel{\to}{n}}}
H^i(X_y, F_{n, y})=0
$$
where $F_{n, y}=F_n|_{X_y}$. If the point $y$ lies in $\bar{C}\backslash C$, or the fibre 
$X_y$ is not irreducible, then what will happen? Is the direct limit  
still zero? Theorem 2 is our answer.

 \begin{theorem} Under the condition of theorem 1, 
  for every point $y$ in $\bar{C}$, 
we have 
$$ \lim_{{\stackrel{\to}{n}}}
H^i(X_y, F_{n, y})=0.
$$
\end{theorem}
$Proof$. If $y$ is contained in $C$  and the fibre $X_y$ in $X$ is irreducible, we
are done. First let $y\in \bar{C}\backslash C$, 
$E=f^{-1}(y),$
we have a short exact sequence  
$$ 0\longrightarrow 
 {\mathcal{O}}_X(-E)
\longrightarrow 
 {\mathcal{O}}_X
\longrightarrow 
 {\mathcal{O}}_E
\longrightarrow 
0.
$$
Tensoring with $F_n$, we have 
$$ 0\longrightarrow 
 F_n(-E)
\longrightarrow 
 F_n
\longrightarrow 
 F_n|_E
\longrightarrow 
0.
$$
If 
$$
\lim_{{\stackrel{\to}{n}}}
H^i(X, F_n(-E))=0,
$$
since we have
$$
\lim_{{\stackrel{\to}{n}}}
H^i(X, F_n)=0, 
$$
then writing the long exact sequence, we get our claim. 

For any fixed $n$, there is a suitable $l$, such that the map 
$$\alpha _2: \quad     H^i(X, F_n)\longrightarrow H^i(X, F_{n+l-1})
$$
 is zero. For this $n$ and $l$, we have a map
$$H^i(X, F_n(-E))
{\stackrel{\alpha}{\longrightarrow}}
 H^i(X, F_{n+l}(-E)). 
$$
$E$ is  component of $D$ (may not be prime) so
the map $\alpha$ can be factored through three maps as follows  
$$H^i(X, F_n(-E))
{\stackrel{\alpha _1}{\longrightarrow}}
 H^i(X, F_n)
{\stackrel{\alpha _2}{\longrightarrow}}
H^i(X, F_{n+l-1})
{\stackrel{\alpha _3}{\longrightarrow}}
H^i(X, F_{n+l}(-E)).
 $$
Since $\alpha _2=0$, and $\alpha=\alpha_3\circ \alpha_2 \circ \alpha_1$, 
we have  $\alpha =0$, i.e, the direct limit we want is zero. 
The map $\alpha_3$
is the natural map corresponding to 
the map 
$F_{n+l-1}\rightarrow 
F_{n+l-1}\otimes {\mathcal{O}}(D-E)$.

If $y$ is a point in $C$ and the fibre $X_y$ is not irreducible, write 
$X_y=X_y'+D'$ where $D'$ is a divisor  contained in $D$, 
$X_y'$ intersects $Y$ with a surface $S$, 
and $X_y'\cap D$ is a curve, 
then $S$ is the fibre of $f|_Y$ over
$y$  in $Y$. $S$ may not be irreducible, however, 
$X_y'\backslash S$ is a divisor on  $X_y'$. By theorem 1 and lemma 4, for every 
$i>0,$
$$
\lim_{{\stackrel{\to}{n}}}
H^i(X'_y, F_n|_{X_y'})=0.
$$
From the short exact sequence
$$ 0\longrightarrow 
 F_n(-X_y')
\longrightarrow 
 F_n
\longrightarrow 
 F_n|_{X_y'}
\longrightarrow 
0
$$
we have for $i+1=2,3,$
$$
\lim_{{\stackrel{\to}{n}}}
H^{i+1}(X, F_n(-X_y'))=0.
$$
Similar to the above argument about $H^i(X, F_n(-E)),$
 we can see that 
$$
\lim_{{\stackrel{\to}{n}}}
H^i(X, F_n(-D'))=0.
$$
For $i>0,$ Consider the map 
$$H^{i+1}(X, F_n(-X_y'-D'))
{\stackrel{\beta}{\longrightarrow}}
 H^{i+1}(X, F_{n+l}(-X_y'-D')). 
$$
As before, it can be factored through three maps as follows
$$H^{i+1}(X, F_n(-X_y'-D'))
{\stackrel{\beta _1}{\longrightarrow}}
 H^{i+1}(X, F_n(-X_y'))
{\stackrel{\beta _2}{\longrightarrow}}
H^{i+1}(X, F_{n+l-1}(-X_y'))
{\stackrel{\beta _3}{\longrightarrow}}
H^{i+1}(X, F_{n+l}(-X_y'-D')).
 $$
For every fixed $n$, we can choose $l$ such that the map $\beta _2$
is zero. Since $\beta=\beta_3\circ \beta_2\circ \beta_1$, $\beta=0$, i.e., for
$i+1=2, 3,$
$$
\lim_{{\stackrel{\to}{n}}}
H^{i+1}(X, F_n(-X_y'-D'))=0.
$$ 
Again from the exact sequence
$$ 0\longrightarrow 
 F_n(-X_y'-D')
\longrightarrow 
 F_n
\longrightarrow 
 F_n|_{X_y'+D'}=F_n|_{X_y}
\longrightarrow 
0
$$
we  get  our claim.
\begin{flushright}
 Q.E.D. 
\end{flushright}

From theorem 1 and theorem 2 we know how the global vanishing cohomology 
controls the local (fibre) cohomology. How 
does the fibre behavior influence the global behavior? 
If for every fibre $S$ of $f|_Y$ in $Y$ over $y\in C$ satisfies 
$H^i(S, \Omega^j_Y|_S)=0$, then  is  $H^i(Y, \Omega^j_Y)=0$? 
We will see that the second and third cohomology vanish but 
the first cohomology is a mystery. To guarantee its vanishing, 
we have to add some mild condition.

To see how the local fibre behavior influences the global behavior, the higher direct images  $R^if_*F_n$ are the link. 
They are coherent for all $i\geq 0$ by Grauert's theorem. Since $f$ is flat 
over $\bar{C}$,  $h^i(X_y, F_{n,y})$=dim$_{\Bbb{C}}H^i(X_y, F_{n, y})$ 
is upper semi-continuous function on  $\bar{C}$. Since 
$H^4(X_y, F_{n, y})=H^3(X_y, F_{n, y})=0$, by [Mu], corollary 3, 
$R^3f_*F_n=0$.
This guarantees  $H^3(Y, \Omega^j_Y)=0$ for every $j$
 and for every point $y\in \bar{C},$
$$R^2f_*F_n{\otimes}
{\Bbb{C}}(y) \cong H^2(X_y, F_{n, y}).
$$ 
If we only consider the closed points $y$ on 
$\bar{C}$,
${\Bbb{C}}(y)=\Bbb{C}$, 
 we have ([U2], Theorem 1.4)
$$(R^2f_*F_n)_y{\otimes}
{\Bbb{C}} \cong H^2(X_y, F_{n, y}),
$$
where   
$(R^2f_*F_n)_y$ is the stalk at $y$ and the tensor product is over ${\mathcal{O}}_{\bar{C},y}$.
So  every stalk satisfies
$$\lim_{{\stackrel{\to}{n}}}
(R^2f_*F_n)_y/{\mathcal{P}}(R^2f_*F_n)_y
 =\lim_{{\stackrel{\to}{n}}}
H^2(X_y, F_{n, y})=0 
$$ 
for every closed point $y$, where ${\mathcal{P}}$ is the maximal ideal of 
${\mathcal{O}}_y$.
This means that for every fixed $n$ and fixed $y$, there is an 
$l$ such that the map
$$\phi:\quad (R^2f_*F_n)_y/{\mathcal{P}}(R^2f_*F_n)_y
\longrightarrow
(R^2f_*F_{n+l})_y/{\mathcal{P}}(R^2f_*F_{n+l})_y
$$ 
is zero. Choose an affine  open neighborhood $U$ of $y$ in $\bar{C}$ such that 
 $R^2f_*F_n|_U=\tilde{M}$,  $R^2f_*F_{n+l}|_U=\tilde{N}$, where $M$ and $N$ are finitely generated modules over $A={\mathcal{O}}(U).$ For every maximal ideal $\mathcal{P}$ of ${\mathcal{O}}(U)$,  we have commutative diagram 
\[
  \begin{array}{ccc}
    M                            &
     {\stackrel{\psi}{\longrightarrow}} &
     N                               \\
    \Big\downarrow\vcenter{%
        \rlap{$\scriptstyle{\pi_1}$}}              &  &
    \Big\downarrow\vcenter{%
       \rlap{$\scriptstyle{\pi_2}$}}      \\
  M/{\mathcal{P}}M              &  {\stackrel{\phi}{\longrightarrow}}  &
  N/{\mathcal{P}}N.   
\end{array}
\]
We can prove  $\psi(M)=0$ if $D_y=D\cap X_y$ is a curve on the fibre $X_y$
for every $y\in C$.  
Therefore
$$
\lim_{{\stackrel{\to}{n}}}
R^2f_*F_n|_C =0 
$$
 which means $H^2(Y, \Omega^j_Y)=0$ by remark 5. 
 In fact, we can get stronger result. From the exact sequence
$$0\longrightarrow 
{\mathcal{O}}(nD)
\longrightarrow
{\mathcal{O}}((n+1)D) 
\longrightarrow
{\mathcal{O}}_D((n+1)D) 
\longrightarrow
0, 
$$
tensoring with $F$ then with ${\mathcal{O}}_{X_y}$, we have 
$$0\longrightarrow 
F_{n,y}
\longrightarrow
F_{n+1, y} 
\longrightarrow 
F_{n+1, y}|_D 
\longrightarrow 
0.
$$
If  $D_y$ is a curve, then  $H^2(F_{n+1, y}|_D)=H^2(F_{n+1, D_y})=0$ 
for every $n$. So the map $H^2(F_{n, y})\rightarrow H^2(F_{n+1, y})$
is surjective for every $n$. But by theorem 2, for suitable $l$, 
the map $H^2(F_{n, y})\rightarrow H^2(F_{n+l, y})$
is zero. Thus there is an $n(y)$ depending on $y$ such that for every 
$n\geq n(y)$,  $H^2(F_{n, y})=0$.  Now fix some $y_0$ in $\bar{C}$ such that 
$H^2(F_{n, y_0})=0$ for every $n\geq n(y_0)$ and there is an open neighborhood
 $U_0$ of $y_0$ in $\bar{C}$ such that $R^2f_*F_{n(y_0)}$  is locally free  
on   $U_0$. Then $H^2(F_{n(y_0), y})=0$ for every $y$ in  $U_0$. So   
$H^2(F_{n, y})=0$ for every $y$ in  $U_0$ and every $n\geq n(y_0)$. 
Let $C-U_0=\{y_1, y_2, ..., y_m\}$, choose 
$n_0=max (n(y_0), n(y_1), ..., n(y_m))$, 
then  
$H^2(X_y, F_{n, y})=0$
for every  $y\in C$ and every $n\geq n_0$.
By  upper semi-continuity theorem,  
$(R^2f_*F_n)_y/{\mathcal{P}}(R^2f_*F_n)_y=0$
for all points $y$ in $C$. By  Nakayama's lemma,   $R^2f_*F_n|_C=0$.

The sheaf $R^1f_*F_n$ is not so nice. For any fixed $n$, 
there is an open set $U_n$ in  
 $\bar{C}$, such that it is locally free on  $U_n$. Let 
$U_n=\bar{C}\backslash A_n$, where $A_n$ is closed in  $\bar{C}$, i.e., 
it consists only finitely many points of  $\bar{C}$.  Since any complete 
metric space is a Baire space (in complex topology, 
every countable intersection of dense 
open sets in  $\bar{C}$ is dense in  $\bar{C}$
  [B2], Chapter 9), 
$B=\bar{C}\backslash \cup A_n=\cap U_n$ is a dense (but we do not know if $B$
is open)
subset of $\bar{C}$ in complex topology. Hence for every point $y$ on $B$,
 all stalks    
$(R^1f_*F_n)_y$ are locally free. Write $B$ as a union of connected subsets 
$B_m$, $B=\cup B_m$, then there is one  $B_m$, such that $B_m$ is dense in 
$\bar{C}$ and connected in complex topology. So we may assume that $B$ 
is connected.  
Again by upper-semicontinuity theorem, for every point $y$ in $C$
 and every $n\geq n_0$, since  $R^2f_*F_n|_C=0$, we have [Mu] 
$$(R^1f_*F_n)_y{\otimes}
{\Bbb{C}} \cong H^1(X_y, F_{n, y}).
$$
For any  $m$, $h^1(X_y, F_{m, y})$ is constant on $B$ since 
$R^1f_*F_m$ is locally free at every point $y$ on $B$
and $B$ is connected. So for the above $n$ and for all points  
$y$ in  $B$, there is $l$ such that  the map 
$$
H^1(X_y, F_{n, y})\longrightarrow H^1(X_y, F_{n+l, y})
$$
is zero.  Moreover, for every point $y$ in $C$ and sufficiently large 
$n$, we have the following 
commutative diagram 

\[
  \begin{array}{ccc}
    R^1f_*F_n{\otimes}
{\Bbb{C}}(y)                            &
     {\stackrel{\approx}{\longrightarrow}} &
     H^1(X_y, F_{n, y})                               \\
    \Big\downarrow\vcenter{%
        \rlap{$\scriptstyle{\alpha}$}}              &  &
    \Big\downarrow\vcenter{%
       \rlap{$\scriptstyle{\beta}$}}      \\
 R^1f_*F_{n+l}{\otimes}
{\Bbb{C}}(y)             &  {\stackrel{\approx}{\longrightarrow}}  &
 H^1(X_y, F_{n+l, y}). 
\end{array}
\]
The map  $\beta$ is zero for every $y\in B$, so as before, the map
$$\alpha :\quad
(R^1f_*F_n)_y/{\mathcal{P}}(R^1f_*F_n)_y
\longrightarrow 
(R^1f_*F_{n+l})_y/{\mathcal{P}}(R^1f_*F_{n+1})_y
$$
is zero for all  points $y$ in $B$.
By the local freeness, this says  on $B$,    
$$\lim_{{\stackrel{\to}{n}}}
R^1f_*F_n|_B =0.  
$$ 
To see this, 
fix a point $y_0$ in $B$, for any sufficiently large $n$
and for the above $l$, choose an affine open set $V$ containing 
$y_0$ such that both $R^1f_*F_n$ and 
 $R^1f_*F_{n+l}$ are locally free on $V$. So there are two 
 positive integers $m_1$ and $m_2$ such that 
$R^1f_*F_n(V)= {\mathcal{O}}(V)^{m_1}$ 
and 
$R^1f_*F_{n+l}(V)= {\mathcal{O}}(V)^{m_2}$. Now for infinitely many 
maximal ideal ${\mathcal{P}}$,  we have commutative 
diagram 
\[
  \begin{array}{ccc}
    {\mathcal{O}}(V)^{m_1}                            &
     {\stackrel{\psi}{\longrightarrow}} &
   {\mathcal{O}}(V)^{m_2}                                \\
    \Big\downarrow\vcenter{%
        \rlap{$\scriptstyle{\pi_1}$}}              &  &
    \Big\downarrow\vcenter{%
       \rlap{$\scriptstyle{\pi_2}$}}      \\
  {\mathcal{O}}(V)^{m_1} /{\mathcal{P}}{\mathcal{O}}(V)^{m_1}              &  
  {\stackrel{\phi}{\longrightarrow}}  &
  {\mathcal{O}}(V)^{m_2}/{\mathcal{P}}{\mathcal{O}}(V)^{m_2}.   
\end{array}
\]
Since  $\psi({\mathcal{O}}(V)^{m_1})
\subset \cap {\mathcal{P}}{\mathcal{O}}(V)^{m_2}=0$, 
where $\mathcal{P}$ runs over infinitely many maximal ideals of 
${\mathcal{O}}(V)$, we have $\psi({\mathcal{O}}(V)^{m_1})=0$. This proves 
$$\lim_{{\stackrel{\to}{n}}}
R^1f_*F_n|_B =0.  
$$

Since the direct limit of $R^1f_*F_n$  
 is quasi-coherent, 
its support is locally closed.  Now
 $B$ is dense and connected in complex topology, there exists an affine open set $U$ 
(we come back to Zariski topology!) in $\bar{C}$
such that on $U$, the direct limit  
$$\lim_{{\stackrel{\to}{n}}}
R^1f_*F_n|_U =0.  
$$
By the following theorem 4 and its remark, we  have proved

\begin{theorem} If for every point $y$ in $\bar{C}$ and for every $i>0$, 
$$ \lim_{{\stackrel{\to}{n}}}
H^i(X_y, F_{n, y})=0,
$$
and for every point $y\in C$,  $D_y=X_y\cap D$ is a curve on the fibre $X_y$,
then   $ R^2f_*F_n|_C =0$
for $n\geq n_0$  and  
$$    \lim_{{\stackrel{\to}{n}}}R^1f_*F_n|_U =0, \quad   
 {\mbox{for a suitable}}\quad U.
$$
 So $H^3(Y, \Omega^j_Y)=
H^2(Y, \Omega^j_Y)=H^1(V, \Omega^j_Y|_V)=0$ for every $j$, 
where $V=f^{-1}(U)\cap Y$. 
\end{theorem}

We have seen  that $R^3f_*F_n=0$ for every $n$
 is determined by the dimension of fibres.  Now we explain why  for $i=1, 2$, 
 if 
$$\lim_{{\stackrel{\to}{n}}}
R^if_*F_n =0,  
$$
then  $H^2(Y, \Omega^j_Y)=H^1(Y, \Omega^j_Y)=0.$
For any point $y\in \bar{C},$ choose an affine open set $U$ containing 
$y$, let $G$ denote the direct limit of $F_n$, then we have long exact sequence of local cohomology
$$
H^1_Z(f^{-1}(U), G)\longrightarrow 
H^1(f^{-1}(U), G)
\longrightarrow
H^1(f^{-1}(U-\{y\}), G)
\longrightarrow   
 H^2_Z(f^{-1}(U), G)
$$
$$
\longrightarrow    H^2(f^{-1}(U), G)
\longrightarrow    H^2(f^{-1}(U-\{y\}), G)
\longrightarrow     H^3_Z(f^{-1}(U), G)
\longrightarrow     0,
$$
where  $Z=f^{-1}(y).$
Since direct limit commutes with cohomology [H1], for every $i>0$, 
$$H^i(f^{-1}(U), G)=
\lim_{{\stackrel{\to}{n}}}
H^i(f^{-1}(U), F_n)
=\lim_{{\stackrel{\to}{n}}}
R^if_*F_n(U)=0 
$$
and 
$$H^i(f^{-1}(U-\{y\}), G)=
\lim_{{\stackrel{\to}{n}}}
H^i(f^{-1}(U-\{y\}), F_n)
=\lim_{{\stackrel{\to}{n}}}
R^if_*F_n(U-\{y\})=0, 
$$
we have $ H^2_Z(f^{-1}(U), G)=  H^3_Z(f^{-1}(U), G)=0.$ Let $V=\bar{C}-\{y\},$ from
$$\longrightarrow 
H^i_Z(X, G)\longrightarrow 
H^i(X, G)
\longrightarrow
H^i(f^{-1}(V), G)
\longrightarrow 
$$
and for $i=2, 3,$
$H^i_Z(X, G)= H^i_Z(U, G)=0.$  By Lemma 4,  we have 
$$ H^2(X, G)=
\lim_{{\stackrel{\to}{n}}}
H^2(X, F_n)=0 \Longrightarrow H^2(Y, \Omega^j_Y)=0,
$$
and
$$ H^3(X, G)=
\lim_{{\stackrel{\to}{n}}}
H^3(X, F_n)=0 \Longrightarrow H^3(Y, \Omega^j_Y)=0.
$$
Now look at $H^1(X, G)$, since
$$\lim_{{\stackrel{\to}{n}}}
R^1f_*F_n(C)=
\lim_{{\stackrel{\to}{n}}}
H^1(f^{-1}(C), F_n)=0
$$
and $Y=f^{-1}(C)\cap Y 
{\stackrel{\psi}{\hookrightarrow}}f^{-1}(C)$ 
is affine morphism ($D$ is locally defined by one equation), by Grothendieck,
[G], page 100, finally we get 
$$0=\lim_{{\stackrel{\to}{n}}}
H^1(f^{-1}(C), F_n)
=H^1(f^{-1}(C), \psi_*(F_n|_Y))
=H^1(Y, \Omega^j_Y).
$$ 
 
\begin{theorem}  If 
$$\lim_{{\stackrel{\to}{n}}}
R^1f_*F_n =\lim_{{\stackrel{\to}{n}}}
R^2f_*F_n =0,
$$
or for every $y\in \bar{C},$ if $D_y=X_y\cap D$ is a curve on the fibre $X_y$
and 
$$\lim_{{\stackrel{\to}{n}}}
R^1f_*F_n
=\lim_{{\stackrel{\to}{n}}}
H^2(X_y, F_{n, y})=0,
$$
then $H^i(Y, \Omega^j_Y)=0$ for every $i>0$ and $j\geq 0$.
\end{theorem}

{\bf{Remark 4}}
 If $R^1f_*F_n$ are locally free, then 
$R^0f_*F_n$ and $R^2f_*F_n$ are locally free by [Mu], page 50, corollary and
from theorem 3 and 4, we know that vanishing direct 
limit of $H^i(X_y, F_{n, y})$ guarantee vanishing of 
Hodge cohomology of $Y$. So it is almost true that the local 
vanishing Hodge cohomology on every fibre $S$ over $C$ guarantee 
the global vanishing of $Y$.  The local freeness of $R^1f_*F_n$  
also tells us that 
$H^i(Y, \Omega^j_Y)=0$ and 
$$\lim_{{\stackrel{\to}{n}}}
R^1f_*F_n =\lim_{{\stackrel{\to}{n}}}
R^2f_*F_n =0
$$
are equivalent.

{\bf{Remark 5}} By lemma 4, the theorem is true if the above 
assumptions hold on $C$, i.e.,  
$$\lim_{{\stackrel{\to}{n}}}
R^1f_*F_n|_C =\lim_{{\stackrel{\to}{n}}}
R^2f_*F_n|_C =0,
$$
or for every $y\in  C,$
$$\lim_{{\stackrel{\to}{n}}}
R^1f_*F_n|_C
=\lim_{{\stackrel{\to}{n}}}
H^2(X_y, F_{n, y})=0.
$$

Now let us consider the affineness of $Y$. Under the 
basic assumption (BA), if $Y$ is affine, then every 
fibre $S$ of $f|_Y$ over $C$ is affine in proposition since it is closed 
in $Y$.  Conversely, if every fibre is 
affine in theorem 1, is  $Y$  affine?
In surface case, it is true. Let us state it precisely. 
If we have a surjective 
morphism from a smooth surface $S$ with $H^i(S, \Omega^j_S)=0$ for every $i>0,
j\geq 0$ to an affine curve $C$, then $S$ must be affine. 
If not, there are 
nonconstant regular functions on $S$ lifted from regular 
functions on $C$. But by Lemma 1.8 [Ku], we know there 
is no  such function on $S$. How about 
the case of 
threefolds? We can give some answer.  
By Serre's affineness criterion, $S$ is affine if and only 
if for all coherent sheaves of ideals
 ${\mathcal{I}}_S$    on $S$,     $H^i(S, {\mathcal{I}}_S)=0$
 for all $i>0,$ or if and only if for  all 
 coherent sheaves ${\mathcal{F}}_S$ on $S$,  $H^i(S, {\mathcal{F}}_S)=0.$
Since the proof of Theorem 2,  theorem 3, and theorem 4 also
works for   coherent sheaves,  we have

 \begin{theorem} (1) In the diagram of proposition, if $Y$ is affine, 
 then every fibre $S$ over $C$ is affine and for every point $y\in \bar{C}$,
every $i>0$ and every 
coherent  sheaf $\mathcal{F}$ on $X$,  
$$\lim_{{\stackrel{\to}{n}}}
H^i(X_y,  {\mathcal{F}}_{n, y})=0,
\quad \quad
\lim_{{\stackrel{\to}{n}}}
R^2f_*{\mathcal{F}}_n =0,
$$
and there is an affine open set $U$ in $\bar{C}$, an integer $n_0$ such that 
for  every $m\geq n_0$, 
$$ \lim_{{\stackrel{\to}{n}}}
R^1f_*{\mathcal{F}}_n|_U=0, \quad \quad
R^2f_*{\mathcal{F}}_m|_U=0,
 $$
where  ${\mathcal{F}}_n=
{\mathcal{F}}\otimes {\mathcal{O}}(nD)$,  ${\mathcal{F}}_{n, y}= 
{\mathcal{F}}_n|_{X_y}.$ 

(2) 
 Conversely, if 
$$\lim_{{\stackrel{\to}{n}}}
R^1f_*{\mathcal{F}}_n|_C =\lim_{{\stackrel{\to}{n}}}
R^2f_*{\mathcal{F}}_n|_C =0,
$$
or for every $y\in C,$ if $D_y$ is a curve and
$$\lim_{{\stackrel{\to}{n}}}
R^1f_*{\mathcal{F}}_n|_C
=\lim_{{\stackrel{\to}{n}}}
H^2(X_y, {\mathcal{F}}_{n, y})=0,
$$
then $Y$ is affine.
\end{theorem}

{\bf{Remark 6}} If $Y$ is affine then it is Stein. 
Theorem 5(2) is also a sufficient condition of Steinness. 
\begin{theorem} If $H^i(Y, \Omega^j_Y)=0$ for every 
$i>0$ and $j\geq 0$, and the $D$-dimension of $X$ is not zero, 
then $Y$ is affine if and only if 
for every coherent sheaf $F$ on $X$,
$$  h^1(X,\lim_{{\stackrel{\to}{n}}}
F\otimes {\mathcal{O}}(nD)) <\infty.
$$
\end{theorem}
$Proof$. By the assumption,  we know that $Y$ contains 
no complete curves. By [GH], 
Proposition 3, we are done.    

\begin{flushright}
 Q.E.D. 
\end{flushright}

For the $D$-dimension  and Kodaira dimension of $X$, we have

\begin{theorem}  If  
$H^i(Y, \Omega^j_Y)=0$, the $D$-dimension of $X$ is not zero, 
and there is a
 smooth fibre $X_{y_0}$ of $f$ over $y_0\in \bar{C}$ such 
 that $S_0=X_{y_0}|_Y$ is not affine, then the Kodaira dimension of $X$
is $-\infty$  and  the $D$-dimension of $X$ is 1.  Generally, we have 
$$  \kappa (\bar{C})+\kappa (X_{y_0}) \leq  
\kappa (X)\leq \kappa (X_{y_0})+ 1.
$$ 
In particular,  if the genus of $\bar{C}$: $g(\bar{C})\geq 2,$ then 
$\kappa (X)=\kappa (X_{y_0})+1.$
\end{theorem}
$Proof$. In the surface case, if $S_0$ is not affine and satisfies 
the same vanishing condition, then the Kodaira dimension of 
its completion  $X_{y_0}$ is  $-\infty$
and the $D$-dimension  is 0 by [I3], [Ku], [Mi].
And if $S_0$ is not affine, then $X_{y_0}$ is birational to
 either the special ruled surface of case (2) or special rational
  ruled surface of case (3) in the first page with $S$ fixed. 
By deformation theorems of Iitaka (I4), (I5), there is an affine open set $U$
in  $\bar{C}$, such that every fibre $X_y$ of $f$ over 
$y\in U$ is of the same 
type.
 By Theorem 5.11 and Theorem 6.12 of Ueno, [U2], we have
$$\kappa (X)\leq \kappa (X_{y_0})+ 1=-\infty.
$$
Combining with upper semicontinuity theorem, if for general fibre $X_y$ over 
$y\in \bar{C}$,  $\kappa (D|_{X_{y_0}}, X_{y_0})=\kappa (D|_{X_{y}}, X_{y})$,
then
$$\kappa (D, X)\leq \kappa (D|_{X_{y_0}}, X_{y_0})+ 1.$$
Consider $\kappa (D|_{X_{y_0}}, X_{y_0})$, if the divisor 
$D_{y_0}=D|_{X_{y_0}}$ on $X_{y_0}$ is a special divisor as [Ku], i.e., 
it has no exceptional divisor of the first type and is  
 a generator of the kernel of the intersection form,
then $H^0(  {\mathcal{O}}_{X_0}(nD_{y_0}))=\Bbb{C}$,  
for every nonnegative integer $n$ (this says  
$H^0({\mathcal{O}}_{X_y}(nD_{y}))=\Bbb{C}$ for every $n$ and general $y$), 
hence $\kappa (D|_{X_{y}}, X_{y})=\kappa (D|_{X_{y_0}}, X_{y_0})=0$,
 $\kappa (D, X)=1.$ But we can not guarantee  that  $D_{y_0}$ 
 is such a special divisor.  By  [I3], properties (1), (2), 
 page 11-12, let $D_1,..., D_r$ be prime components of $D$, 
 for all integers $p_1$,$..., p_r>0,$  we have, 
$$\kappa (D_1+...+D_r,  X)= \kappa (p_1D_1+...+p_rD_r, X);
$$
and if $g: W\rightarrow V$ is surjective morphism,
where $W$ and $V$ are smooth projective varieties, 
 $E$ is an effective divisor on $W$ such that codim$(g(E))\geq 2,$ then
$$\kappa (g^*(D')+E,  W)=\kappa (D', V) 
$$
where $D'$ is a divisor on $V$, $g^*(D')=\sum D_i'$ 
is the reduced transform of $D'$, 
where
$D_i'$ are  irreducible components. 
By [U2], Lemma 5.3, page 51-52, if every fibre 
in the above map $g$ is connected, 
then we have $\Bbb{C}$-linear isomorphism 
$$  H^0(V, {\mathcal{O}}_V(D')) \cong 
H^0(W, {\mathcal{O}}_W(g^*D')).
$$ 
From these properties we get the same $D$-dimension.
In fact, on the fibre $X_y$, the $D_y$-dimension does not depend on the 
support of the divisor $D_y$, i.e., $D_y$  may 
contain exceptional curves of the first kind. It also does not depend 
on the coefficients of the prime divisors of $D_y$.   
 In any case, the open part $S=X_y\cap Y$ is fixed. Therefore
  $\kappa (D|_{X_{y_0}}, X_{y_0})=0.$   Hence
 $\kappa (D, X)=1$ if $S_{y_0}$ is not affine. 

The left  cases  follow  from [Ka2], [V].

\begin{flushright}
 Q.E.D. 
\end{flushright}

{\bf{Remark 7}} \quad In the above proof, if the fibre $X_{y_0}=X_0$ is smooth,
then $S_0=X_0\cap Y$ is smooth and satisfies 
$H^i(S_0, \Omega^j_{S_0})=0$. If $S_0$ is not affine, then it is 
fixed by Mohan Kumar's classification, i.e., it is either type (2) surface 
in the 
first page
or  type (3) surface. Here the boundary $D_0=X_0-S_0=D|_{X_0}$ may 
not be the special 
divisor $D_0'$ in [Ku]. But by the above argument, 
$\kappa (D_0, X_{0})=\kappa (D_0', X_{0})=0.$  It might happen
that there is no any global divisor $D$ on $X$ such that when restricted to 
the fibre $X_0$, it is the divisor $D_0'$. Fortunately, we do not need 
the existence of such a special  divisor $D$. In fact, we have $D$ first 
then consider 
its restriction on the fibre.

  In surface case, if  $H^i(S, \Omega^j_S)=0$ for all 
  $i>0$ and $j\geq 0,$  $S$ is not affine, then the 
  Kodaira dimension of its completion
is unique and the $D$-dimension is also unique. 
In threefold case, is it still  true?   Using the 
same notation as in Theorem 7, if $S_0$ is affine, 
then we can choose $U$ such that 
every fibre $X_y$  of $f$ over $U$ has constant Kodaira dimension. 
If $g(\bar{C})\geq 2$, then 
$\kappa (X)=\kappa (X_y)+1.$ But under this condition, is $Y$ affine? 
This is 
equal to the question that  if every (or general) fibre S of 
$f|_Y$ over $C$ is affine, 
is $Y$ affine? Usually it is not true even in surface case 
without the restriction of Hodge cohomology.

   Consider the logarithmic Kodaira dimension 
$\bar{\kappa}(Y)=\kappa(K_X+D, X)$. 
For general fibre $X_y$, $D|_{X_y}=D_y$ is a divisor
 on $X_y$ with normal crossings [I1]. The 
 logarithmic Kodaira dimension does not depend on the 
 embedding if the boundary is a divisor with normal crossings [Mi]. 
By theorem 3 [I3], $\bar{\kappa}(Y)\leq 
\bar{\kappa}(S)+1$, where $S$ is general fibre. 
If $S$ is smooth but not affine, then 
  $\bar{\kappa}(S)=-\infty$,
therefore $\bar{\kappa}(Y)=-\infty$. 
Generally, 
 if Iitaka's $\bar{C}_n$ 
conjecture is true [I3], i.e., 
$\bar{\kappa}(Y)\geq \bar{\kappa}(S)+\bar{\kappa}(C),$ then 
$$       \bar{\kappa}(S)+\bar{\kappa}(C) \leq \bar{\kappa}(Y)
\leq \bar{\kappa}(S)+1.
$$ 
In particular, if $\bar{\kappa}(C)=1$, i.e., 
the genus $g(\bar{C})\geq 2,$ then 
$\bar{\kappa}(Y)=\bar{\kappa}(S)+1.$
In 1978, Kawamata [Ka3] proved this conjecture if the fibre dimension is 1.

Finally we give an example.

{\bf{Example }} Let $S$ be a surface with 
$H^i(S, \Omega^j_S)=0$ for all $i>0$ and $j\geq 0,$ not affine, let $C$ 
be any affine curve, then $Y=S\times C$ satisfies $H^i(Y, \Omega^j_Y)=0$ 
for all $i>0$ and $j\geq 0$ by K$\mbox{\"{u}}$nneth formula (see [Hi] or [SW]). 
The $D-$dimension of 
$X=\bar{S}\times \bar{C}$ is 1 and the Kodaira dimension is $-\infty$ by 
theorem 7. 
Its logarithmic Kodaira dimension is  also $-\infty$. 
   Again by  K$\mbox{\"{u}}$nneth formula, 
$q(X)=h^1({\mathcal{O}}_X)=g(\bar{C}).$
  This example says that $q(X)$ can be any 
nonnegative  integer. So if we choose different $C$ with different genus, the 
corresponding $X$ are not isomorphic since they have different $q$.

We are constructing  nonproduct examples and will submit it later.  
\\   
\\

Department of Mathematics, Washington University, St. Louis, MO 63130

E-mail address: zhj@@math.wustl.edu

\end{document}